\documentclass{article}
\usepackage{amssymb}
\usepackage{amsmath}
\newtheorem{mytheorem}{Theorem}
\title{The Zipf law for random texts with unequal probabilities of occurrence of letters and the Pascal pyramid}
\author{V.V.\,Bochkarev and E.Yu.\,Lerner}
\date{}
\begin{document}

\maketitle

\begin{abstract}
We model the generation of words with independent 
unequal probabilities of occurrence of letters. 
We prove that the probability $p(r)$ of occurrence of words of rank~$r$ has a power asymptotics. 
As distinct from the paper published earlier by B.~Conrad and M.~Mitzenmacher, 
we give a brief proof by elementary methods and obtain an explicit formula for the exponent of the power law.
\end{abstract}

\textbf{Keywords:}\ {Zipf law, monkey model, order statistics, power laws, Pascal pyramid, recursive sequences, functional equations.}

\medskip

As is known, in English texts the word ``the'' occurs most often, the next one (with respect to the occurrence) is the word ``of'', therefore, each word-form~$w$ in a text is associated with its {\it rank} $r(w)$, i.e., the number in the frequency list. The {\it frequency} $f(w)$ of the word $w$ in a text is defined as the ratio of the number of occurrences of the word~$w$ to the length of this text. According to the Zipf law (established in the first half of the last century), the product of $r(w)$ and $f(w)$ approximately equals a constant value; for an English text it equals $0.1$. At the present time, owing to Google Labs, available are results of the recognition of 4\% of all books ever published~\cite{google}. According to our numerical tests, the OLS line constructed by logarithmically transformed (for convenience, we use the decimal logarithm) values of $r$~and~$f$ of one hundred English words used most often in 2000, takes the form $\lg f=-1.05182 - 1.00026\lg r$. The fact that the modern data (of the indicated sample) agree with the Zipf law is the starting point of our work.

The Zipf law was interpreted in many ways; see~\cite{maslov} for a brief review of relevant papers published in Russia and the inference of this law from the general properties of semiotic systems proposed by V.P.~Maslov. In this paper we consider a more traditional explanation of the Zipf law, namely, the ``monkey'' model. We assume that a ``monkey'' independently with equal probabilities types any of 26~English letters or does the space with the probability~$p_0$~\cite{Mand}, \cite{Miller},~\cite{Li}. We understand a word as a sequence of letters between two spaces. Evidently, occurrences of all words of the same length have equal probabilities, and their ranks run in succession. We denote by~$p(r)$ the probability of obtaining a word of the rank~$r$. Evidently, $\exists c_1,c_2$: $c_1<\ln p(r)-\alpha\ln r<c_2$, where in the case $p_0=1/27$ we have $\alpha=\ln 27/\ln 26$. Really, various letters in a text have various frequencies of occurrence; in Russian they are approximately determined by the law established by S.~Gusein-Zade~\cite{Gus1}, \cite{Gus2} (probabilities are proportional to average values of order statistics of the exponential distribution).

Not long time ago B.~Conrad and M.~Mitzenmacher~\cite{Mit} have generalized the ``monkey'' models for the case when letters have different probabilities of occurrence. They proved an inequality analogous to~(\ref{BochkarevLernerTwo}) below by using the Tauber theorems for generating functions. In conclusion of the mentioned papers the authors write ``It would, of course, be pleasant to have a proof of the power law behavior in the case of unequal probabilities that avoids some of this technical machinery'' and discuss the possible ways to obtain such a result. In this paper we propose a short proof obtained in another way, namely, by making use of properties of the Pascal pyramid. Note that in the case of a Markovian dependence of probabilities of letters the power law does not necessarily holds. We study the power law (as well as the exponential one) for Markovian chains (and obtain explicit formulas for the corresponding parameters) in a separate paper.

\begin{mytheorem}
\it Let probabilities of letters equal $p_1,\ldots,p_n$, $n>1$ ($\sum_{i=1}^n p_i=1-p_0$), and let $\gamma$ be a root of the equation $\sum_{i=1}^n p_i^\gamma=1$. Then
\begin{equation}
\label{BochkarevLernerTwo}
\exists c_1,c_2:\qquad c_1< \ln p(r)-\ln r/\gamma < c_2.
\end{equation}
\end{mytheorem}

{\bf Proof.}
Without loss of generality, for convenience, we introduce an empty word; we assume that its rank equals one, while other words have greater ranks. All words that contain $k_1$ letters of the 1st kind, $\ldots$, $k_n$ words of the $n$th kind have one and the same probability $\mathop{\rm Pr}(k)=p_1^{k_1}\ldots p_n^{k_n} p_0$, and their ranks run in succession. The number of such words is defined by the multinomial coefficient~$M(k_1,\ldots,k_n)$,
$M(k_1,\ldots,k_n)=\frac{(k_1+\ldots+k_n)!}{k_1!\ldots k_n!}$. Let us fix some probability~$f$ ($f\in (0,1]$) and denote by $Q(f)$ the rank of the last word~$w$, whose probability is not less than~$f$ in the sorted (in the nonincreasing order of probabilities) in the list of all words. Thus, for example, $Q(p_0)=1$ (we take into account only the empty word), $Q(p' p_0)=2$, where $p'=\max\{p_1,\ldots,p_n\}$ (here we assume that the maximum is unique),~etc. Evidently, the function $Q$ is nonincreasing, piecewise constant (it takes on only positive integer values), and tending to infinity as $f\to 0$.

The rank of a word~$w$ (i.e., the number of words at the beginning of the mentioned list up to $w$ inclusive) equals $M(k_1,\ldots,k_n)$, where the sum is taken over all $n$-tuples such that $\mathop{\rm Pr}(k)\ge f$.

Let $\tilde Q(x)=Q(p_0 e^{-x})$, $x\ge 0$. Evidently, $\tilde Q(x)$ is
\begin{equation}
\label{BochkarevLernerThree}
\sum\limits_{k\ge 0: L_1 k_1+\ldots L_n k_n\le x}
M(k_1,\ldots,k_n),
\end{equation}
where $L_i=-\ln p_i$, $i=1,\ldots,n$. We have to prove that the value $\ln\tilde Q(x)-\gamma x$ is bounded. This is a property of a multidimensional generalization of the Pascal pyramid. Note that we have $\sum_{i=1}^n e^{-L_i}=1-p_0$. Introduce $L'_i=\gamma L_i$ and assume that $x'=\gamma x$. We have $\sum_{i=1}^n e^{-L'_i}=1$, and correlation~(\ref{BochkarevLernerThree}) takes the form 
$$\sum _{k\ge 0: L'_1 k_1+\ldots L'_n k_n\le x'}
M(k_1,\ldots,k_n).
$$
We have to prove that by subtracting $x'$ from the logarithm of this sum we obtain a bounded value. Thus, the general case is reduced to the case when~$\gamma=1$.

Let us now immediately prove the boundedness of the difference, assuming that 
\begin{equation}
\label{BochkarevLernerFour}
\sum_{i=1}^n e^{-L_i}=1.
\end{equation}
The function $\tilde Q(x)$ is defined for $x\ge 0$. For convenience, we extend it with zero for~${x<0}$. Analogously we extend the Pascal pyramid with zero for negative $k_i$. Then any number in the Pascal pyramid, except $M(0,\ldots,0)$, is the sum of successive numbers such that one of their indices is less by one. This leads to the functional equation $\tilde Q(x)=\tilde Q(x-L_1)+\ldots+\tilde Q(x-L_n)+Q_0(x)$, where $Q_0$ is the Heaviside step function (it equals zero for negative values of the argument and does one for its nonnegative values).

For $x\ge L'=\max\{ L_1,\ldots,L_n\}$ we obtain the recurrent correlation
\begin{equation}
\label{BochkarevLernerOne}
\tilde Q'(x)=\tilde Q'(x-L_1)+\ldots+\tilde
Q'(x-L_n),
\end{equation}
where $\tilde Q'(x)=\tilde Q(x)+1/(n-1)$. If all fractions $L_i/L_j$, $i,j=1,\ldots,n$, are rational, then we obtain a recursive sequence, and one can prove the theorem with the help of known explicit formulas for such sequences. In a general case, instead of nuclear formulas we use simple inequalities.

Multiplying correlation~(\ref{BochkarevLernerFour}) by $e^x$, we obtain that the exponential function satisfies Eq.~(\ref{BochkarevLernerOne}). Since the function $\tilde Q'(x)$ is piecewise constant, the product $q(x)=\tilde Q'(x) e^{-x}$ is piecewise continuous, and hence there exist positive $c_1$ and $c_2$ such that $c_1< q(x)< c_2$ for all $0\le x\le L'$. In the recurrent correlation~(\ref{BochkarevLernerOne}) we replace addends in the right-hand side with their lower (upper) bounds and thus widen the interval, where the inequality $c_1< q(x)< c_2$ is fulfilled, up to $x\le L'+L''$, where $L''=\min\{ L_1,\ldots,L_n\}$. Repeating this procedure several times, in a finite number of steps we prove the inequality for any arbitrarily large~$x$. By applying the logarithmic transform to the inequality we conclude that $\ln \tilde Q'(x)-x$ is bounded, and then so is the difference $\ln \tilde Q(x)-x$, which was to be proved. 
\medskip

Thus, we have proved the power order of the asymptotics and obtained an explicit formula for it. In the case of the function~$Q(x)$ the parameter of the power law is defined by the solution~$\gamma$ of the equation $\sum_{i=1}^n p_i^\gamma=1$. But if we consider the probability~$p(r)$ of the occurrence the $r$th word in the frequency list, then we obtain the power order of the asymptotics with the exponent $1/\gamma$ (i.e., greater than one). Note that the verification of the power law with respect to the frequency of occurrence of one thousand (rather than one hundred) most often cited words in several European languages included in the database of Google Labs (French, German, Russian, Spain, and several variants of English) shows that the exponent for each language essentially differs from one and is close to the value calculated by our formula.

We are grateful to V.D.~Solovyev for useful discussions.


\begin{thebibliography}{9}
\bibitem{google} 
J.B.\,Michel, Y.K.\,Shen, A.P.\,Aiden, A.~Veres, M.K.\,Gray, J.P.\,Pickett, D.~Hoiberg, D.~Clancy, P.~Norvig, J.~Orwant, S.~Pinker, M.A.\,Nowak, and E.L.\,Aiden, ``Quantitative Analysis of Culture Using Millions of Digitized Books,'' Science, \textbf{331} (6014), 176--182 (2011),
http://www.librarian.net/wp-content/uploads/science-googlelabs.pdf.

\bibitem{maslov}
V.P.\,Maslov and T.V.\,Maslova, ``On Zipf's Law and Rank Distributions in Linguistics and Semiotics,'' Matem. Zametki \textbf{80} (5),  718--732 (2006).

\bibitem{Mand}
B.~Mandelbrot, ``An Informational Theory of the Statistical Structure of Languages,'' in \textquotedblleft Communication Theory, W. Jackson, Ed: \textquotedblright\ (Betterworth, 1953), pp.~486--502.

\bibitem{Miller}
G.A.\,Miller, ``Some Effects of Intermittent Silence,'' 
Amer. J. Psychology \textbf{70}, 311--314 (1957).

\bibitem{Li}
W.Li, ``Random Texts Exhibit Zipf's-Law-Like Word Frequency Distribution,'' IEEE Transactions on Information Theory {\bf 38} 1842--1845 (1992).

\bibitem{Gus1}
S.M.\,Gusein-Zade, ``Frequency Distribution of Letters in the Russian Language,'' Probl. Peredachi Inform. \textbf{24} (4), 102-107 (1989).

\bibitem{Gus2}
S.M.\,Gusein-Zade, 
``On the Frequency of Meeting of Key Words and on Other Ranked Series,''
Nauchno-Tekhnicheskaya Informatsiya, Ser.~2, N~1, 28--31 (1987).

\bibitem{Mit}
B.~Conrad and M.~Mitzenmacher, 
``Power Laws for Monkeys Typing Randomly: The Case of Unequal Probabilities,'' IEEE Transac. \textbf{50}, 1403--1414 (2004),
http://www.eecs.harvard.edu/\~{}michaelm/postscripts/toit2004a.pdf.
\end{thebibliography}
\end{document}